\documentclass[11pt,a4paper,twoside]{amsart}
\usepackage{amsfonts, amssymb,indentfirst}

\newtheorem{prop}{Proposition}[section]
\newtheorem{theorem}[prop]{Theorem}

\newtheorem{defi}[prop]{Definition}

\newtheorem{remark}[prop]{Remark}

\newcommand{\cqd}{\hfill$\Box$}
\renewcommand{\hom}[0]{\operatorname{Hom}}

\newcommand{\quot}[0]{\operatorname{Quot}}

\renewcommand{\dim}[0]{\operatorname{dim}}
\newcommand{\codim}[0]{\operatorname{codim}}

\newcommand{\rk}[0]{\operatorname{rk}}

\title[INTERSECTION NUMBERS ON THE COMPACT VARIETY OF RATIONAL RULED SURFACES]
{INTERSECTION NUMBERS ON THE COMPACT VARIETY OF RATIONAL RULED
SURFACES}
\author[Cristina Mart\'inez]{ Cristina Mart{\'\i}nez Ram\'irez}
\address{Max Planck Institute for Mathematics, Vivatgasse 7, Bonn, 53111.}
 \email{cmartine@mpim-bonn.mpg.de}

\begin{document}
\maketitle

\begin{abstract}

We describe a natural action on the Quot scheme, $R_{d}$
compactifying the space of degree $d$ maps from $\mathbb{P}^{1}$
to the Grassmannian of lines. We identify the fixed points
components for this action and the weights of the normal bundle of
these components. We compute the degree of this Quot scheme under
the generalized Pl\"ucker embedding by applying Atiyah-Bott
localization formula. 

\end{abstract}


\section{Introduction.}

The Quot scheme is a fine moduli space equipped with an universal
element. It has been used many times as a smooth compactification
of the space of morphisms of a fixed degree from a curve $C$ to a
Grassmannian, \cite{Ber1}. It is known that when $C$ is of genus
0, the Quot scheme is irreducible and smooth. Recalling the
notation of \cite{Mar1}, we denote by $R_{d}$ the Quot scheme
parametrizing quotients of rank 2 and degree $d$ of a trivial
vector bundle $\mathcal{O}^{4}_{\mathbb{P}^{1}}$, and by
$R^{0}_{d}$ the open set of morphisms. S.A. Str\o mme in
\cite{Str} computes the Betti numbers of $R_{d}$ and gives a
description with ge\-ne\-ra\-tors and relations of its cohomology
ring. In particular, he gives a
basis for the Picard group $A^{1}(R_{d})$. 

Here we compute the degree of the generalized Pl\"ucker embedding
of the Quot scheme $R_{d}$. It is called generalized Pl\"ucker
embedding because in some sense can be considered a generalization
of the Pl\"ucker embedding given by the hyperplane class of the
corresponding Grassmannian. In our case, we are considering the
Grassmannian of lines in $\mathbb{P}^{3}$, which we will denote as
$G(2,4)$, but similar computations can be obtained for a general
Grassmannian $G(k,n)$ of $k-$planes in $\mathbb{C}^{n}$. In
particular for $n-k=1$, the Quot scheme is a projective space and
therefore intersection theory here is well understood. The method
that we use here, involves the geometry of the Quot scheme and the
Bott residue formula. 

The Bott residue formula expresses the degree of certain
zero-cycles classes on a smooth complete variety with an action of
an algebraic torus in terms of local contributions supported on
the components of the fixed point set. We describe a natural
action on the Quot scheme. We identify the fixed points components
for this action and the weights of the normal bundle of these
components. We study the equivariant codimension one cohomology
group of $R_{d}$ for this action.

\subsubsection*{Notation.}
 We work over the field of complex numbers $\mathbb{C}$. Let $X$ be a projective
 homogeneous variety, then $A_{d}X$ and $A^{d}X$ can be
 taken to be the Chow homology and cohomology groups of $X$. We identify $A^{d}X$ with $A_{n-d}X$
 by the Poincar\'e duality isomorphism. We use cup pro\-duct
$\cup$ for the product in $A^{*}X$. The moduli space
$\overline{M}_{0,n}(X,\beta)$ parametrizes marked stable maps from
genus 0 curves to $X$ in the cohomology class $\beta \in
H^{2}(X,\mathbb{Z})$. Let $\gamma_{i}$ be cycles on $X$, if $\sum
\codim(\gamma_{i})=\dim(\overline{M}_{0,n}(X,d))$, {\it the
Gromov-Witten invariant}

\noindent $I_{0,n,d}(\gamma_{1},\ldots,\gamma_{n})$ is defined as
the top degree class:
\begin{equation}
I_{0,n,d}(\gamma_{1},\ldots,\gamma_{n})=\int_{[\overline{M}_{0,n}(X,d)]^{vir}}
\pi_{1}^{*}(\gamma_{1})\cup \ldots \cup \pi_{n}^{*}(\gamma_{n}),
\end{equation}
where $[\overline{M}_{0,n}(X,d)]^{vir}$ denotes the virtual
fundamental class of $\overline{M}_{0,n}(X,d)$.

\section{The degree of the generalized Pl\"ucker embedding and the Vafa-Intriligator formula.}
When we fix the degree, $d$, and the rank, 2, of a locally free
sheaf $E$ on $\mathbb{P}^{1}$, we are fixing its Hilbert
polynomial,

\begin{equation}\label{eq10}
P(t)=\chi(E(t))=d+2t+2.
\end{equation}

The moduli, $\quot_{d}(\mathbb{P}^{1},P(t))$, of quotients with
fixed Hilbert polynomial $P(t)$, is a fine moduli space which we
will denote as $R_{d}$. We observe the quotient
$\mathcal{O}^{4}_{\mathbb{P}^{1}}\rightarrow E\rightarrow 0$,
determines a point $q\in \quot(\mathbb{P}^{1},P(t))$ and a
morphism $f_{q}:\mathbb{P}^{1}\rightarrow G(2,4)$, by the
universal property of the Grassmannian. By definition, there is an
universal quotient

\begin{equation}\label{eq11}
\mathcal{O}^{4}_{\quot\times\mathbb{P}^{1}}\rightarrow
\mathcal{E}_{\quot\times \mathbb{P}^{1}}.
\end{equation}
 S.A. Str\o mme in \cite{Str} gives a basis for the
Picard group $A^{1}(R_{d})$, formed by the divisors:

$$\alpha=c_{1}(\pi_{1*}(\mathcal{E}\otimes
\pi_{2}^{*}\mathcal{O}_{\mathbb{P}^{1}}(d))-c_{1}(\pi_{1*}(\mathcal{E}\otimes
\pi_{2}^{*}\mathcal{O}_{\mathbb{P}^{1}}(d-1))),$$
$$\beta=c_{1}(\pi_{1*}(\mathcal{E}\otimes
\pi_{2}^{*}\mathcal{O}_{\mathbb{P}^{1}}(d-1)),$$

\noindent where $\mathcal{E}$ is the universal quotient over
$R_{d}\times \mathbb{P}^{1}$, and $\pi_{1}$, $\pi_{2}$ are the
projection maps over the first and second factors respectively.
Let $e:R^{0}_{d}\times\mathbb{P}^{1} \rightarrow G(2,4)$ be the
evaluation map and $T_{1}\in H_{6}(G(2,4),\mathbb{Z})$, $T_{a}\in
H_{4}(G(2,4),\mathbb{Z})$ be the classes of an hyperplane and
$a-$plane respectively. The hyperplane class determines the
Pl\"ucker embedding of the Grassmannian $G(2,4)$ as a quadric in
$\mathbb{P}^{5}$ which corresponds to a variety of lines in
$\mathbb{P}^{3}$. The $a-$plane $T_{a}$ corresponds to lines in
$\mathbb{P}^{3}$ contained in a given plane. The following set of
morphisms define Weil divisors on $R^{0}_{d}$:
$$ \label{div1}
D_{i}:=\{\varphi \in R^{0}_{d}|\ \  e(t,\varphi)\cap T_{a_{i}}\neq
\emptyset\},
$$
$$\label{div2} Y:=\{\varphi \in R^{0}_{d}|\
e(t_{i},\varphi)\in T_{1}\ \mbox{for  a  fixed } \  t_{i}\in
\mathbb{P}^{1} \}.
$$
These divisors extend to divisors on $R_{d}$.

Let $P_{d}$ be the degree of $R_{d}$ by the morphism induced by
the polarization given by the divisor $\alpha$, by lemma 3.2 of
\cite{Mar2}, $P_{d}$ is the degree of the top codimensional
cohomology class given by the autointersection
\begin{equation}
P_{d}=\int _{[R_{d}]}[Y]^{4d+4}\cap
[R_{d}]=\int_{[R_{d}]}\alpha^{4d+4}\cap [R_{d}].
\end{equation}

This intersection number is computed in \cite{RRW} via Quantum
Cohomology. Note that in this case, the intersection is
transversal in the Quot scheme compactification, and therefore the
intersection number is the same than the given by integrating over
the space of stable maps, that is the Gromov-Witten invariant
$I_{0,4d+4,d}(T_{1},\stackrel{4d+4}{\ldots},T_{1})$. It can be
obtained too by means of the formulas of Vafa and Intriligator,
\cite{Ber1}. This does not happen when considering intersection
numbers containing $D_{i}$.

\subsubsection*{ Vafa and Intriligator's Formula:} Let
$\zeta$ be a primitive $n$th root of $(-1)^{k}$ and assume that
$0\leqq a_{i}\leqq k$ and $a_{1}+ \ldots +
a_{N}=\dim\,(\overline{M}_{0,4d+4}(G(k,n),d))$. Then:

$$I_{0,4d+4,d}(\gamma_{a_{1}},\stackrel{4d+4}{\ldots},\gamma_{a_{N}})=(-1)^{k\choose
2}n^{-k}\sum_{i_{1}>\ldots>i_{k}}\sigma_{a_{1}}(\zeta^{I})\ldots\sigma_{a_{N}}(\zeta^{I})
\frac{\prod_{j\neq
l}(\zeta^{i_{j}}-\zeta^{i_{l}})}{\prod_{j=1}^{k}\zeta^{(n-1)i_{j}}} $$

\noindent where $\zeta^{I}=(\zeta^{i_{1},\ldots,\zeta^{i_{k}}})$ and
$\sigma_{a_{i}}$ are the elementary symmetric polynomials in $k$ variables.

\vspace{0.5cm} In \cite{Mar1} it is proved that the
top-codimensional classes involving the divisors $D_{i}$ have not
enumerative meaning on the Quot scheme because there is an excess
component of intersection contained in the boundary. Then the
intersection is carried out in the Kontsevich compactification of
stable maps, $\overline{M}_{0,3}(G(2,4),d)$.


\noindent The tool we are going to use for computing these
degrees, is the Bott residue formula. For this purpose, we
consider the equivariant action of a one dimensional torus
$T=\mathbb{C}^{*}$ on the variety $R_{d}$. First we study the
varieties of fixed points under this action and we compute the
equivariant Chern classes in the Chow equivariant rings of the
normal bundles restricted to the varieties of fixed points and the
line bundles corresponding to the divisors we intersect.

\subsubsection*{Bott's residues formula,}(see \cite{cf}).

\noindent Let $X$ be a smooth, complete variety of dimension $n$
and let $E$ be a $T-$equivariant vector bundle over $X$. Then we
have
\begin{equation} \label{Bott}
\int_{X}(p(E)\cap
[X])=\sum_{F\subset\,R^{T}}\pi_{F*}\left(\frac{p^{T}(E_{|_{F}})\cap
[F]_{T}} {c_{d_{F}}(\mathcal{N}_{F/X})}\right). \end{equation}
Here $d_{F}$ denotes the codimension of the component $F$ of fixed
points in $X$ and  $p^{T}(E)$ is the polynomial of degree $n$ in
the line bundles corresponding to the cycles expressing the
product in the equivariant Chow ring of $F$. The numerator will be
a product of polynomials of degree 1 corresponding to the line
bundles restricted to $F$. The denominator will be a polynomial
with the dimension of the normal bundle as degree.

\noindent In our case, we will apply the formula to compute the
intersection
$$\int_{[R_{d}]}\alpha^{4d+4}\cap [R_{d}].$$


\section{Varieties of fixed points under the
$\mathbb{C}^{*}-$ action}
 We consider the diagonal action of the one-dimensional torus acting
 on the variety $R_{d}$.

\vspace{0.5cm} \noindent A point in $R_{d}$ is given by a
quotient: $0\rightarrow N\rightarrow \mathcal{O}^{4}\rightarrow
E\rightarrow 0$ in $\mathbb{P}^{1}$, \ \ where
$\chi(E(t))=2t+2+d$. Let $\overline{w}=(w_{0},w_{1},w_{2},w_{3})$
be a quadruple of integral weights with $w_{0}<w_{1}<w_{2}<w_{3}$.
$\mathbb{C}^{*}$ acts on each point:
\[\forall t\in \mathbb{C}^{*}, \begin{array}{ccc}
\mathbb{C}^{*}\times \mathbb{C}^{4} & \rightarrow &
\mathbb{C}^{4}\\ t\cdot (x_{0},x_{1},x_{2},x_{3})= & &
(t^{w_{0}}x_{0},t^{w_{1}}x_{1},t^{w_{2}}x_{2},t^{w_{3}}x_{3})
\end{array}\]
Let $0\rightarrow \mathcal{N}\rightarrow \mathcal{O}^{4}_{
R_{d}\times \mathbb{P}^{1}}\rightarrow \mathcal{E}\rightarrow 0$
be the universal exact sequence in $R_{d}\times\mathbb{P}^{1}$.
This action induces an isomorphism $\mathcal{O}^{4}\rightarrow
\mathcal{O}^{4}$, such that the weight corresponding to the
trivial sheaf $\mathcal{O}_{R_{d}}$ is $w_{0}+w_{1}+w_{2}+w_{3}$,
since $\pi_{*}\wedge^{4}\mathcal{O}^{4}_{R_{d}}\cong
\mathcal{O}_{R_{d}}$, where $\pi: R_{d}\times
\mathbb{P}^{1}\rightarrow R_{d}$.

\vspace{0.5cm} \noindent Let $E$ be of rank 2 with a nonzero
torsion of degree $d$. We suppose,
\begin{equation} \label{zd}E\cong \mathcal{O}_{Z_{d}}\oplus 0\oplus
\mathcal{O}_{\mathbb{P}^{1}}\oplus \mathcal{O}_{\mathbb{P}^{1}}.
\end{equation}

\noindent The scheme $Hilb^{d}_{\mathbb{P}^{1}}$ will denote the
Hilbert scheme of $d$ points in $\mathbb{P}^{1}$ which is
isomorphic to $\mathbb{P}^{d}$ and $\mathcal{Z}_{d}\subset
\mathbb{P}^{d}\times \mathbb{P}^{1}$ will denote the incidence
variety. The scheme $Hilb^{d}_{\mathbb{P}^{1}}$ parametrizes
points $\mathcal{O}^{4}_{\mathbb{P}^{1}}\rightarrow E\rightarrow
0$, with $E$ as in (\ref{zd}).

\begin{defi}
Let $\underline{P}=(P_{i})_{1\leq i\leq r}$ be a family of
polynomials with rational coefficients. $\underline{P}$ is said to
be a good partition of the polynomial $\underline{P}$, if
$Hilb^{P_{i}}_{\mathbb{P}^{1}}\neq \emptyset$ and
$\underline{P}=\sum_{1\leq i\leq r}P_{i}$.
\end{defi}

\subsubsection*{Examples of good partitions of the Hilbert
polynomial of E.}

\begin{enumerate}
\item $P_{1}(t)=d,\ P_{2}(t)=0,\ P_{3}(t)=t+1,\ P_{4}(t)=t+1$, is
said to be a {\it good partition} of the polynomial
$\underline{P}(t)=2t+2+d$.

\item Another good partition for $d$ odd would be:
\begin{equation}P_{\frac{d-1}{2},\frac{d+1}{2}}(t)=\underbrace{\frac{d+1}{2}}+
\underbrace{\frac{d-1}{2}}+\underbrace{t+1}+\underbrace{t+1}\end{equation}
and for $d$ even:
\begin{equation} \label{d2} P_{\frac{d}{2},\frac{d}{2}}(t)=
\underbrace{t+1}+\underbrace{t+1}+\underbrace{\frac{d}{2}}+
\underbrace{\frac{d}{2}} \end{equation}
\end{enumerate}
$\mathcal{O}^{4}_{\mathbb{P}^{1}}\rightarrow
\mathcal{O}_{Z_{d/2}}\oplus \mathcal{O}_{Z_{d/2}}\oplus
\mathcal{O}_{\mathbb{P}^{1}}\oplus \mathcal{O}_{\mathbb{P}^{1}}
\rightarrow 0$, is a boundary point corresponding to the partition
(\ref{d2}).
$$Hilb_{\mathbb{P}^{1}}^{t+1}\times
Hilb_{\mathbb{P}^{1}}^{t+1}\times
Hilb_{\mathbb{P}^{1}}^{d/2}\times Hilb_{\mathbb{P}^{1}}^{d/2}$$
is the subscheme of the scheme
$Quot(\mathcal{O}^{4}_{\mathbb{P}^{1}},P(t))$ that we will denote
by $(t+1,t+1,\frac{d}{2},\frac{d}{2})$ in order to simplify
notation. The fixed point set of $R_{d}$ under the action of the
torus $T$ is the union of the possible subschemes
$Hilb^{\underline{P}}_{\mathbb{P}^{1}}$ \cite{Bif}. 

\subsection{Fixed points in $R_{d}$ under the
$\mathbb{C}^{*}$-action}

\begin{prop}
The varieties of fixed points in $R_{d}$ under the
$\mathbb{C}^{*}-$action are parametrized by,
\[\mathbb{P}^{d}\]\[\mathbb{P}^{d-1}\times \mathbb{P}^{1}
\] \[\mathbb{P}^{d-2}\times\mathbb{P}^{2}\]\[\vdots\]
\[\left\{\begin{array}{ll}\mathbb{P}^{\frac{d+1}{2}}\times
\mathbb{P}^{\frac{d-1}{2}} & \mbox{if $d$ odd}\\
\mathbb{P}^{\frac{d}{2}}\times \mathbb{P}^{\frac{d}{2}} & \mbox{if
$d$ even.}\end{array}\right.\] There are 12 components of each
type.
\end{prop}
{\it Proof.} Following the work of Bifet \cite{Bif}, we see that
to study the components of fixed points under the
$\mathbb{C}^{*}-$action, is equivalent to study good partitions
for the Hilbert polynomial $\underline{P}(t)=2t+2+d$.
\begin{enumerate} \item Corresponding to the partition,
\[P_{d,0}(t)=d+0+\underbrace{t+1}+\underbrace{t+1},\] we have,
\[Hilb^{d}_{\mathbb{P}^{1}} \times Hilb^{0}_{\mathbb{P}^{1}}\times
Hilb^{t+1}_{\mathbb{P}^{1}}\times Hilb^{t+1}_{\mathbb{P}^{1}}\cong
\mathbb{P}^{d}. \] There are 12 of this kind.

\item Corresponding to the partition,
\[P_{b,a}(t)=b+a+\underbrace{t+1}+\underbrace{t+1},\ \ \ b\geq a>0,\
b+a=d,\]\[Hilb^{b}_{\mathbb{P}^{1}} \times
Hilb^{a}_{\mathbb{P}^{1}}\times Hilb^{t+1}_{\mathbb{P}^{1}}\times
Hilb^{t+1}_{\mathbb{P}^{1}}\cong \mathbb{P}^{b}\times
\mathbb{P}^{a} \] There are $\frac{d}{2}-1$ different components
if $d$ even, and $\frac{d+1}{2}-1$ if $d$ odd. These are
parametrized by:\[\mathbb{P}^{d-1}\times \mathbb{P}^{1}
\] \[\mathbb{P}^{d-2}\times\mathbb{P}^{2}\]\[\vdots\]
\[\left\{\begin{array}{ll}\mathbb{P}^{\frac{d+1}{2}}\times
\mathbb{P}^{\frac{d-1}{2}} & \mbox{if $d$ odd}\\
\mathbb{P}^{\frac{d}{2}}\times \mathbb{P}^{\frac{d}{2}} & \mbox{if
$d$ even.}\end{array}\right.\] For each one there are also ${4
\choose 2}\times 2$ components.
\end{enumerate}
There are $12\cdot \frac{d}{2}=6d$ fixed point components if $d$
even, and $12\cdot \frac{d+1}{2}=6\cdot (d+1)$ if $d$ odd. \cqd


\vspace{1cm} {\bf The Euler characteristic of $R_{d}$} is given by
the formula (see corollary 1.4 of \cite{Str}),
\[\chi(R_{d})=\rk_{\mathbb{Z}}\,A\,(R_{d})={4\choose
2}{d+4-1\choose d}.\] Since the Chow ring of $\mathbb{P}^{d}$ is
$\mathbb{Z}[h]/\langle h^{d}\rangle$ the contribution to the Euler
characteristic of the component of fixed points of the first kind
is $d+1$ cycles. The contribution of the components of the second
kind to the Euler characteristic is $a\cdot b$, since
$A^{*}(\mathbb{P}^{b}\times
\mathbb{P}^{a})=\mathbb{Z}(H,h)/\langle H^{b},h^{a}\rangle$.

\subsubsection*{Example: fixed points in $R_{3}$ under the
$\mathbb{C}^{*}$-action} The Quot scheme $R_{3}$ parametrizes
quotients with Hilbert polynomial,
\[P_{3}(t)=2t+2+3.\] There are two kinds of fixed points varieties in $R_{3}$
corresponding to the two partitions of the polynomial $P_{3}(t).$
\begin{enumerate}
\item \[P_{3,0}(t)=3+0+\underbrace{t+1}+\underbrace{t+1}\]
\[Hilb^{3}_{\mathbb{P}^{1}}\times Hilb^{0}_{\mathbb{P}^{1}}\times
Hilb^{t+1}_{\mathbb{P}^{1}}\times Hilb^{t+1}_{\mathbb{P}^{1}}\cong
\mathbb{P}^{3}.\] There are ${4 \choose 2}\times 2=12$ of this
kind, each one contributes 4 cycles, since the Chow ring of
$\mathbb{P}^{3}$ is $A^{*}(\mathbb{P}^{3})\cong
\mathbb{Z}(h)/\langle h^{4}\rangle$, here $h$ is the class of a
hyperplane in $\mathbb{P}^{3}$. 

\item
\[P_{2,1}(t)=2+1+\underbrace{t+1}+\underbrace{t+1}\]
\[Hilb^{2}_{\mathbb{P}^{1}}\times Hilb^{1}_{\mathbb{P}^{1}}\times
Hilb^{t+1}_{\mathbb{P}^{1}}\times Hilb^{t+1}_{\mathbb{P}^{1}}\cong
\mathbb{P}^{2}\times \mathbb{P}^{1}.\] There are ${4 \choose
2}\times 2=12$ of this kind, and each one contributes 6 cycles in
$A^{*}(\mathbb{P}^{2}\times \mathbb{P}^{1})\cong
\mathbb{Z}(H,h)/\langle H^{3},h^{2}\rangle$, where $H$ is the
class of a hyperplane in $\mathbb{P}^{2}$ and $h$ in
$\mathbb{P}^{1}$.
\end{enumerate}
This number coincides with the Euler characteristic, since:
\[120=\chi\,(R_{3})={4 \choose 2}{3+2(4-2)-1 \choose 3}\]
There are 24 components of fixed points.

\vspace{1cm} \noindent We have seen that the position of the
trivial sheaves and the torsion sheaves is important. It
determines different components up to isomorphism. It also
determines the weights that act.

\subsection{The $\mathbb{C}^{*}$-equivariant Chern classes}
We shall associate to the component
$Hilb^{d}_{\mathbb{P}^{1}}\times Hilb^{0}_{\mathbb{P}^{1}}\times
Hilb^{t+1}_{\mathbb{P}^{1}}\times Hilb^{t+1}_{\mathbb{P}^{1}}\cong
\mathbb{P}^{d}$ the quadruple $(d,0,t+1,t+1)$ to simplify
notation, and $c^{T}_{1}(\alpha|_{(d,0,t+1,t+1)})$,
$c^{T}_{1}(\beta|_{(d,0,t+1,t+1)})$ will be the corresponding
first $\mathbb{C}^{*}$-equivariant Chern classes in the
equivariant Chow ring $A^{T}_{*}(\mathbb{P}^{d})$.
\begin{theorem} The first equivariant Chern classes of the fixed points components are

\noindent  $c^{T}_{1}(\alpha|_{(t+1,t+1,d,0)})=h+w_{0}+w_{1}, \ \
\ c^{T}_{1}(\beta|_{(t+1,t+1,d,0)})=d\,w_{0}+d\,w_{1}+d\,w_{2},$

$c^{T}_{1}(\alpha|_{(t+1,d,t+1,0)})=h+w_{0}+w_{2}, \ \ \
c^{T}_{1}(\beta|_{(t+1,d,t+1,0)})=d\,w_{0}+d\,w_{1}+d\,w_{2},$

$c^{T}_{1}(\alpha|_{(d,0,t+1,t+1)})=h+w_{2}+w_{3}, \ \ \
c^{T}_{1}(\beta|_{(d,0,t+1,t+1)})=d\,w_{0}+d\,w_{2}+d\,w_{3},$

$c^{T}_{1}(\alpha|_{(t+1,t+1,0,d)})=h+w_{0}+w_{1}, \ \ \
c^{T}_{1}(\beta|_{(t+1,t+1,0,d)})=d\,w_{0}+d\,w_{1}+d\,w_{3},$

$c^{T}_{1}(\alpha|_{(d,t+1,0,t+1)})=h+w_{1}+w_{3}, \ \ \
c^{T}_{1}(\beta|_{(d,t+1,0,t+1)})=d\,w_{0}+d\,w_{1}+d\,w_{3},$

$c^{T}_{1}(\alpha|_{(d,t+1,t+1,0)})=h+w_{1}+w_{2}, \ \ \
c^{T}_{1}(\beta|_{(d,t+1,t+1,0)})=d\,w_{0}+d\,d\,w_{1}+d\,w_{2},$

$c^{T}_{1}(\alpha|_{(t+1,t+1,0,d)})=h+w_{0}+w_{1}, \ \ \
c^{T}_{1}(\beta|_{(t+1,t+1,0,d)})=d\,w_{0}+d\,w_{1}+d\,w_{3},$

$c^{T}_{1}(\alpha|_{(t+1,0,t+1,d)})=h+w_{0}+w_{2}, \ \ \
c^{T}_{1}(\beta|_{(t+1,0,t+1,d)})=d\,w_{0}+d\,w_{2}+d\,w_{3},$

$c^{T}_{1}(\alpha|_{(0,d,t+1,t+1)})=h+w_{2}+w_{3}, \ \ \
c^{T}_{1}(\beta|_{(0,d,t+1,t+1)})=d\,w_{1}+d\,w_{2}+d\,w_{3},$

$c^{T}_{1}(\alpha|_{(t+1,0,d,t+1)})=h+w_{0}+w_{3}, \ \ \
c^{T}_{1}(\beta|_{(t+1,0,d,t+1)})=d\,w_{0}+d\,w_{2}+d\,w_{3},$

$c^{T}_{1}(\alpha|_{(0,t+1,d,t+1)})=h+w_{1}+w_{3}, \ \ \
c^{T}_{1}(\beta|_{(0,t+1,d,t+1)})=d\,w_{1}+d\,w_{2}+d\,w_{3},$

$c^{T}_{1}(\alpha|_{(0,t+1,t+1,d)})=h+w_{1}+w_{2}, \ \ \
c^{T}_{1}(\beta|_{(0,t+1,t+1,d)})=d\,w_{1}+d\,w_{2}+d\,w_{3},$


$c^{T}_{1}(\alpha|_{(t+1,t+1,a,b)})=H+h+w_{0}+w_{1}, \ \ \
c^{T}_{1}(\beta|_{(t+1,t+1,a,b)})=aH+bh+dw_{0}+dw_{1}+aw_{2}+bw_{3},$

$c^{T}_{1}(\alpha|_{(t+1,a,t+1,b)})=H+h+w_{0}+w_{2}, \ \ \
c^{T}_{1}(\beta|_{(t+1,a,t+1,b)})=aH+bh+dw_{0}+dw_{1}+aw_{1}+bw_{3},$

$c^{T}_{1}(\alpha|_{(a,b,t+1,t+1)})=H+h+w_{2}+w_{3}, \ \ \
c^{T}_{1}(\beta|_{(t+1,a,t+1,b)})=aH+bh+dw_{0}+aw_{1}+dw_{1}+bw_{3},$

$c^{T}_{1}(\alpha|_{(t+1,a,b,t+1)})=H+h+w_{0}+w_{3}, \ \ \
c^{T}_{1}(\beta|_{(t+1,a,b,t+1)})=aH+bh+dw_{0}+aw_{1}+bw_{1}+dw_{3},$

$c^{T}_{1}(\alpha|_{(a,t+1,b,t+1)})=H+h+w_{0}+w_{3}, \ \ \
c^{T}_{1}(\beta|_{(a,t+1,b,t+1)})=aH+bh+aw_{0}+dw_{1}+bw_{2}+dw_{3},$

$c^{T}_{1}(\alpha|_{(a,t+1,t+1,b)})=H+h+w_{1}+w_{2}, \ \ \
c^{T}_{1}(\beta|_{(a,t+1,t+1,b)})=aH+bh+aw_{0}+dw_{1}+dw_{2}+bw_{3}.$

$c^{T}_{1}(\alpha|_{(t+1,t+1,b,a)})=H+h+w_{0}+w_{1}, \ \ \
c^{T}_{1}(\beta|_{(t+1,t+1,b,a)})=aH+bh+dw_{0}+dw_{1}+bw_{2}+aw_{3}.$

$c^{T}_{1}(\alpha|_{(t+1,a,t+1,b)})=H+h+w_{0}+w_{1}, \ \ \
c^{T}_{1}(\beta|_{(t+1,a,t+1,b)})=aH+bh+dw_{0}+aw_{1}+dw_{2}+bw_{3}.$

$c^{T}_{1}(\alpha|_{(a,b,t+1,t+1)})=H+h+w_{1}+w_{2}, \ \ \
c^{T}_{1}(\beta|_{(a,b,t+1,t+1)})=aH+bh+aw_{0}+bw_{1}+dw_{2}+dw_{3}.$

$c^{T}_{1}(\alpha|_{(t+1,a,b,t+1)})=H+h+w_{0}+w_{3}, \ \ \
c^{T}_{1}(\beta|_{(t+1,a,b,t+1)})=aH+bh+dw_{0}+aw_{1}+bw_{2}+dw_{3}.$

$c^{T}_{1}(\alpha|_{(a,t+1,b,t+1)})=H+h+w_{1}+w_{3}, \ \ \
c^{T}_{1}(\beta|_{(a,t+1,b,t+1)})=aH+bh+aw_{0}+dw_{1}+bw_{2}+dw_{3}.$

$c^{T}_{1}(\alpha|_{(a,t+1,t+1,b)})=H+h+w_{1}+w_{2}, \ \ \
c^{T}_{1}(\beta|_{(a,t+1,t+1,b)})=aH+bh+aw_{0}+dw_{1}+dw_{2}+bw_{3}.$

\end{theorem}

{\it Proof.}
 We consider the universal quotient in $R_{d}\times
\mathbb{P}^{1}$ restricted to the fixed point component
$(d,0,t+1,t+1)$. It is enough to consider one component of fixed
points by the symmetry of the computations.
$$ 0\rightarrow
\mathcal{O}_{\mathbb{P}^{d}\times \mathbb{P}^{1}}(-1,-d)\oplus
\mathcal{O}_{\mathbb{P}^{d}\times \mathbb{P}^{1}} \rightarrow
\mathcal{O}^{4}_{\mathbb{P}^{d}\times \mathbb{P}^{1}}\rightarrow
\mathcal{O}_{\mathcal{Z}_{d}}\oplus 0\oplus
\mathcal{O}_{\mathbb{P}^{d}\times \mathbb{P}^{1}}\oplus
\mathcal{O}_{\mathbb{P}^{d}\times \mathbb{P}^{1}}\rightarrow 0.
$$
$$
\begin{picture}(100,90)
\put(0,70){$\mathbb{P}^{3}\times \mathbb{P}^{1}\supset
\mathcal{Z}_{d}$} \put(15,65){\vector(3,-2){30}}
\put(60,65){\vector(0,-1){25}} \put(65,50){$\pi_{1}$}
\put(60,25){$\mathbb{P}^{d}$}
\end{picture}
$$

\noindent Let $h$ denote the positive generator of the projective
space $\mathbb{P}^{d}$.

\noindent Since $\mathcal{Z}_{d}$ is $\pi_{1}-$flat, the
restriction of the coherent sheaf $B_{d-1}$ to the fixed point
component $(d,0,t+1,t+1)$ is well defined by,

\[
B_{d-1}|_{(d,0,t+1,t+1)}=\pi_{1*}[(\mathcal{O}_{\mathcal{Z}_{d}}\oplus
0 \oplus \mathcal{O}_{\mathbb{P}^{d}\times \mathbb{P}^{1}}\oplus
 \mathcal{O}_{\mathbb{P}^{d}\times
\mathbb{P}^{1}})\otimes
\pi_{2}^{*}\mathcal{O}_{\mathbb{P}^{1}}(d-1)]=\] \[
 \pi_{1}^{*}\mathcal{O}_{\mathcal{Z}_{d}}\otimes\pi_{2}^{*}\mathcal{O}_{\mathbb{P}^{1}}(d-1)
 \oplus \pi_{1*}\pi_{2}^{*}\mathcal{O}_{\mathbb{P}^{1}}(d-1)
 \oplus \pi_{1*}\pi_{2}^{*}\mathcal{O}_{\mathbb{P}^{1}}(d-1), \  \mbox{where}
\]
\[\pi_{1*}\pi_{2}^{*}\mathcal{O}_{\mathbb{P}^{1}}(d-1)=\mathcal{O}^{d}_{\mathbb{P}^{d}}.\]

\noindent We consider the exact sequence defining the sheaf
$\mathcal{O}_{\mathcal{Z}_{d}}$:
\begin{equation} 0\rightarrow \mathcal{O}_{\mathbb{P}^{d}\times \mathbb{P}^{1}}(-1,-d)
\rightarrow \mathcal{O}_{\mathbb{P}^{d}\times
\mathbb{P}^{1}}\rightarrow
\mathcal{O}_{\mathcal{Z}_{d}}\rightarrow 0,
\end{equation}
\noindent and we tensorize with the line bundle
$\pi_{2}^{*}\mathcal{O}_{\mathbb{P}^{1}}(d-1)$ and take the long
exact sequence associated to the pushforward $\pi_{1*}$:

\begin{equation} 0\rightarrow
\pi_{1*}(\mathcal{O}_{\mathbb{P}^{d}\times
\mathbb{P}^{1}}(-1,-d)\otimes
\pi_{2}^{*}\mathcal{O}_{\mathbb{P}^{1}}(d-1)) \rightarrow
\pi_{1*}(\mathcal{O}\otimes\pi_{2}^{*}\mathcal{O}_{\mathbb{P}^{1}}(d-1))\rightarrow
\end{equation}\[\rightarrow
\pi_{1*}(\mathcal{O}_{\mathcal{Z}_{d}}\otimes\pi_{2}^{*}\mathcal{O}_{\mathbb{P}^{1}}(d-1))\rightarrow
0,
\]

\noindent The vanishing of
$R^{1}\pi_{1*}(\mathcal{O}_{\mathbb{P}^{d}\times
\mathbb{P}^{1}}(-1,-1))$ implies that

$\pi_{1*}(\mathcal{O}_{\mathcal{Z}_{d}}\otimes
\pi_{2}^{*}\mathcal{O}_{\mathbb{P}^{1}}(d-1))=\mathcal{O}^{d}_{\mathbb{P}^{d}}$.
Therefore the rank of $B_{d-1}|_{(d,0,t+1,t+1)}$ is $2d+d=3d$ and
the restriction of $\beta_{d}$ to the component $(d,0,t+1,t+1)$ is
$\mathcal{O}_{\mathbb{P}^{d}}$. Since
$\bigwedge^{d}\mathcal{O}^{d}_{\mathbb{P}^{d}}=\mathcal{O}_{\mathbb{P}^{d}}$
the corresponding weight is $dw_{0}+dw_{2}+dw_{3}$.

\noindent The first equivariant Chern class
$\alpha_{d}|_{(d,0,t+1,t+1)}$ is defined by

\noindent
$\alpha_{d}|_{(d,0,t+1,t+1)}=c_{1}(B_{d}|_{(d,0,t+1,t+1)})-c_{1}(B_{d-1}|_{(d,0,t+1,t+1)})$,
with
\begin{align*}B_{d}|_{(d,0,t+1,t+1)} =\pi_{1*}[(\mathcal{O}_{\mathcal{Z}_{d}}\otimes
0\oplus \mathcal{O}_{\mathbb{P}^{d}\times
\mathbb{P}^{1}}\oplus\mathcal{O}_{\mathbb{P}^{d}\times
\mathbb{P}^{1}})\otimes
\pi^{*}_{2}\mathcal{O}_{\mathbb{P}^{1}}(d)] & \\ =
\pi_{1*}\mathcal{O}_{\mathcal{Z}_{d}}\otimes
\pi_{2}^{*}\mathcal{O}_{\mathbb{P}^{1}}(d)\oplus
\pi_{1*}\pi_{2}^{*}\mathcal{O}_{\mathbb{P}^{1}}(d)\oplus
\pi_{1*}\pi_{2}^{*}\mathcal{O}_{\mathbb{P}^{1}}(d).\end{align*} We
have the following exact sequence

\begin{equation} 0\rightarrow
\pi_{1*}(\mathcal{O}_{\mathbb{P}^{d}\times
\mathbb{P}^{1}}(-1,0))\rightarrow
\pi_{1*}(\mathcal{O}\otimes\pi_{2}^{*}\mathcal{O}_{\mathbb{P}^{1}}(d))\rightarrow
\pi_{1*}(\mathcal{O}_{\mathcal{Z}_{d}}\otimes\pi_{2}^{*}\mathcal{O}_{\mathbb{P}^{1}}(d))
\rightarrow 0.
\end{equation}
\noindent We have that
$\pi_{1*}\pi_{2}^{*}\mathcal{O}_{\mathbb{P}^{1}}(d)=\mathcal{O}^{d+1}_{\mathbb{P}^{d}}$
and $R^{1}\pi_{1*}(\mathcal{O}_{\mathbb{P}^{d}\times
\mathbb{P}^{1}}(-1,0))=0$. We also have that
$\pi_{1*}(\mathcal{O}_{\mathbb{P}^{d}\times
\mathbb{P}^{1}}(-1,0))= \mathcal{O}_{\mathbb{P}^{d}}(-1)$,
therefore,
\[\pi_{1*}(\mathcal{O}_{\mathcal{Z}_{d}}\otimes
\pi_{2}^{*}\mathcal{O}_{\mathbb{P}^{1}}(d))=\mathcal{O}_{\mathbb{P}^{d}}(1)\oplus
\mathcal{O}^{d-1}_{\mathbb{P}^{d}}.\]The rank of
$B_{d}|_{(d,0,t+1,t+1)}$ is $2(d+1)+d=3d+2$ and the weight
$dw_{0}+(d+1)w_{2}+(d+1)w_{3}$, therefore the restriction of
$\alpha_{d}$ to $(d,0,t+1,t+1)$ is $$\fbox{$h+w_{2}+w_{3}$}$$


We now compute {\bf the restrictions of the divisors $\alpha_{d}$
and $\beta_{d}$} to the components of fixed points of the second
kind. Consider the following incidence variety
$$
\begin{picture}(100,90) \label{p3}
\put(0,80){$\mathbb{P}^{b}\times \mathbb{P}^{a}\times
\mathbb{P}^{1} \supset \mathcal{Z}_{b}\times \mathcal{Z}_{a}$}
\put(95,75) {\vector(0,-1){25}}\put(100,60){$\pi_{12}$}
\put(80,35) {$\mathbb{P}^{b}\times \mathbb{P}^{a}$}
\put(100,30){\vector(1,-1){15}} \put(100,30){\vector(-1,-1){15}}
\put(80,5){$\mathbb{P}^{b}$}\put(110,5){$\mathbb{P}^{a}$}
\put(75,25){$p_{1}$} \put(110,25){$p_{2}$}
\put(90,75){\vector(-3,-2){40}}\put(30,35){$\mathbb{P}^{1}$}
\put(50,65){$\pi_{3}$}\end{picture} $$ and the restriction of the
coherent sheaf $B_{d-1}$ to the fixed point component
$(b,a,t+1,t+1)$,
\[B_{d-1}|_{(b,a,t+1,t+1)}=\pi_{12*}[(\mathcal{O}_{\mathcal{Z}_{b}}
\oplus \mathcal{O}_{\mathcal{Z}_{a}}\oplus
\mathcal{O}_{\mathbb{P}^{b}\times \mathbb{P}^{a}\times
\mathbb{P}^{1}}\oplus\mathcal{O}_{\mathbb{P}^{b}\times
\mathbb{P}^{a}\times \mathbb{P}^{1}})\otimes
\pi_{3}^{*}\mathcal{O}_{\mathbb{P}^{1}}(d-1)]\  \]\[ \mbox{where}
\ b+a=d,\ b\geq a>0, \ \mbox{and}\ \ \] \[
\pi_{12*}(\mathcal{O}_{\mathbb{P}^{b}\times \mathbb{P}^{a}\times
\mathbb{P}^{1}}(0,-1,-a)\otimes
\pi_{3}^{*}\mathcal{O}_{\mathbb{P}^{1}}(d-1))=\pi_{12*}\mathcal{O}_{\mathbb{P}^{b}\times
\mathbb{P}^{a}\times \mathbb{P}^{1}}(0,-1,b-1).\] Since
$R^{1}\pi_{12*}\mathcal{O}_{\mathbb{P}^{b}\times
\mathbb{P}^{a}\times \mathbb{P}^{1}}(0,-1,b-1)=0$, the following
exact sequence stands: \[0\rightarrow \mathcal{O}(0,-1)\otimes
\mathcal{O}^{b}\rightarrow \mathcal{O}^{d}\rightarrow
\pi_{12*}(\mathcal{O}_{\mathcal{Z}_{a}}\otimes
\pi_{3}^{*}\mathcal{O}_{\mathbb{P}^{1}}(d-1))\rightarrow 0\
\rm{over}\ \mathit{\mathbb{P}^{b}\times \mathbb{P}^{a}}.\] Let
$h,H$ denote the positive generators of the projective spaces
$\mathbb{P}^{a}$ and $\mathbb{P}^{b}$, respectively. The first
Chern class of the bundle
$\pi_{12*}(\mathcal{O}_{\mathcal{Z}_{a}}\otimes
\pi_{3}^{*}\mathcal{O}_{\mathbb{P}^{1}}(d-1))$ is $bh$. \noindent
For computing the first Chern class of the subbundle
$\pi_{12*}(\mathcal{O}_{\mathcal{Z}_{b}}\otimes
\pi_{3}^{*}\mathcal{O}_{\mathbb{P}^{1}}(d-1))$, we see that
\[\pi_{12*}(\mathcal{O}_{\mathbb{P}^{b}\times \mathbb{P}^{a}\times
\mathbb{P}^{1}}(0,-1,-b)\otimes
\pi_{3}^{*}\mathcal{O}_{\mathbb{P}^{1}}(d-1))=\pi_{12*}\mathcal{O}_{\mathbb{P}^{b}\times
\mathbb{P}^{a}\times \mathbb{P}^{1}}(-1,0,a-1), \ \mbox{and},\]
\[R^{1}\pi_{12*}\mathcal{O}_{\mathbb{P}^{b}\times \mathbb{P}^{a}\times
\mathbb{P}^{1}}(-1,0,a-1)=0, \ \mbox{therefore},\]
\[0\rightarrow
\mathcal{O}(-1,0)\otimes \mathcal{O}^{a}\rightarrow
\mathcal{O}^{d}\rightarrow \pi_{12*}(\mathcal{O}_{Z_{b}}\otimes
\pi_{3}^{*}\mathcal{O}_{\mathbb{P}^{1}}(d-1))\rightarrow 0\
\rm{on}\ \mathit{\mathbb{P}^{b}\times \mathbb{P}^{a}}.\] By
symmetry with the previous case, the first Chern class of the
bundle $\pi_{12*}(\mathcal{O}_{\mathcal{Z}_{b}}\otimes
\pi_{3}^{*}\mathcal{O}_{\mathbb{P}^{1}}(d-1))$ is $aH$. It follows
that

\noindent $c_{1}(B_{d-1}|_{(b,a,t+1,t+1)})=aH+bh$, and its weight
is $bw_{0}+aw_{1}+dw_{2}+dw_{3}$. Finally, the restriction of
$\beta_{d}$ to the component $(b,a,t+1,t+1)$, isomorphic to
$\mathbb{P}^{b}\times \mathbb{P}^{a}$, is

$$\fbox{$aH+bh+bw_{0}+aw_{1}+dw_{2}+dw_{3}$}$$

\noindent The restriction of $B_{d}$ to $(b,a,t+1,t+1)$ is given
by
\[B_{d}|_{(b,a,t+1,t+1)}=\pi_{12*}[(\mathcal{O}_{\mathcal{Z}_{b}}\oplus
\mathcal{O}_{\mathcal{Z}_{a}}\oplus
\mathcal{O}\oplus\mathcal{O})\otimes
\pi_{3}^{*}\mathcal{O}_{\mathbb{P}^{1}}(d)].\] We see that
\[\pi_{12*}(\mathcal{O}_{\mathbb{P}^{b}\times \mathbb{P}^{a}\times
\mathbb{P}^{1}}(0,-1,-a)\otimes
\pi_{3}^{*}\mathcal{O}_{\mathbb{P}^{1}}(d))=\pi_{12*}\mathcal{O}_{\mathbb{P}^{b}\times
\mathbb{P}^{a}\times \mathbb{P}^{1}}(0,-1,b)\ \mbox{and}\]
\[R^{1}\pi_{12*}\mathcal{O}_{\mathbb{P}^{b}\times \mathbb{P}^{a}\times
\mathbb{P}^{1}}(0,-1,b)=0.\]

\noindent Therefore we have the exact sequence:\[0\rightarrow
\mathcal{O}(0,-1)\otimes \mathcal{O}^{b+1}\rightarrow
\mathcal{O}^{d+1}\rightarrow
\pi_{12*}(\mathcal{O}_{\mathcal{Z}_{a}}\otimes
\pi_{3}^{*}\mathcal{O}_{\mathbb{P}^{1}}(d))\rightarrow 0\
\rm{over}\ \mathit{\mathbb{P}^{b}\times \mathbb{P}^{a}}. \]

\noindent Again from the fact that
\[ \pi_{12*}(\mathcal{O}_{\mathbb{P}^{b}\times \mathbb{P}^{a}\times
\mathbb{P}^{1}}(-1,0,-b)\otimes
\pi_{3}^{*}\mathcal{O}_{\mathbb{P}^{1}}(d))=\pi_{12*}\mathcal{O}_{\mathbb{P}^{b}\times
\mathbb{P}^{a}\times \mathbb{P}^{1}}(-1,0,a), \  \]
\[\ \mbox{and} \ \ \mbox{that} \ \ R^{1}\pi_{12*}\mathcal{O}_{\mathbb{P}^{b}\times \mathbb{P}^{a}\times
\mathbb{P}^{1}}(-1,0,a)=0,\] it follows that there is an exact
sequence,
\[0\rightarrow
\mathcal{O}(-1,0)\otimes \mathcal{O}^{a+1}\rightarrow
\mathcal{O}^{d+1}\rightarrow
\pi_{12*}(\mathcal{O}_{\mathcal{Z}_{b}}\otimes
\pi_{3}^{*}\mathcal{O}_{\mathbb{P}^{1}}(d))\rightarrow 0\
\rm{over}\ \mathit{\mathbb{P}^{b}\times \mathbb{P}^{a}}.\]
\noindent The first Chern class of $B_{d}|_{(b,a,t+1,t+1)}$ is
$(a+1)H+(b+1)h$ and its weight
$bw_{0}+aw_{1}+(d+1)w_{2}+(d+1)w_{3}$, therefore the restriction
of $\alpha_{d}$ to $(b,a,t+1,t+1)$ is, $$
\fbox{$H+h+w_{2}+w_{3}$}$$
\cqd

\section{The normal bundle.}\label{fnormal}
\begin{theorem}
 The Equivariant Chern classes of the Normal bundle in the
equivariant Chow ring of the fixed points component are, for the
first kind of components:

\begin{align*}
c_{3d+4}^{T}(\mathcal{N}_{\mathbb{P}^{d}/R_{d}})=(h+(w_{2}-w_{0}))^{d+1}\cdot
 (h+(w_{3}-w_{0}))^{d+1}\cdot(w_{0}-w_{1})&\\\cdot
 (h-(w_{0}-w_{1}))^{d-1}(w_{2}-w_{1})\cdot(w_{3}-w_{1}),
\end{align*}
and for the second kind of components:


\noindent
\begin{align*}c_{3d+4}^{T}(\mathcal{N}_{\mathbb{P}^{b}\times
\mathbb{P}^{a}}/R_{d})=(H+w_{2}-w_{0})^{b+1}\cdot
(H+w_{3}-w_{0})^{b+1}&\\ \cdot (h+w_{2}-w_{0})^{a+1}\cdot
(-H+h+w_{0}-w_{1})^{b-a-1}\cdot  (h+w_{0}-w_{1})^{a+1} &\\ \cdot
(H+w_{1}-w_{0})^{b+1}\cdot (H-h+w_{1}-w_{0})^{a-b-1},\end{align*}
for $b-a\geq 0$.

\end{theorem}

{\it Proof.} We consider the normal bundle of the fixed points
component in $R_{d}$.

\noindent We need to compute the weight of the normal bundle for
each component of fixed points, and its equivariant Chern class.
We first study the normal bundle for the components of the first
kind.

\noindent We consider again the universal exact sequence

\[0\rightarrow \mathcal{N}\rightarrow \mathcal{O}^{4}_{R_{d}\times
\mathbb{P}^{1}}\rightarrow \mathcal{E}\rightarrow 0.\] The tangent
space to the variety $R_{d}$ is  $\mathcal{T}_{R_{d}}\cong
\pi_{*}\mathcal{H}om\left(\mathcal{N},\mathcal{E}\right)$ (see
$\S$7.7.1 of \cite{Str}), where $\pi:R_{d}\times
\mathbb{P}^{1}\rightarrow R_{d}$. If we restrict it to a component
of fixed points of the first kind:

\[0\rightarrow \mathcal{Y}_{1}\oplus \mathcal{Y}_{2}\oplus
\mathcal{Y}_{3}\oplus\mathcal{Y}_{4}\rightarrow
\mathcal{O}^{4}_{\mathbb{P}^{d}\times \mathbb{P}^{1}}\rightarrow
\mathcal{E}_{\mathbb{P}^{d}\times \mathbb{P}^{1}}\rightarrow 0,\]

\noindent where $\bigoplus_{i=1}^{4}\mathcal{Y}_{i}$ is the kernel
of the quotient map $\mathcal{O}^{4}_{\mathbb{P}^{d}\times
\mathbb{P}^{1}}\rightarrow \mathcal{E}_{\mathbb{P}^{d}\times
\mathbb{P}^{1}}\rightarrow 0$. The restriction of the normal
bundle to the component of fixed points yields
\[0\rightarrow \mathcal{T}_{\mathbb{P}^{d}}\rightarrow
\mathcal{T}_{R}\rightarrow \mathcal{N}_{\mathbb{P}^{d}}\rightarrow
0,\] \[\mathcal{N}_{\mathbb{P}^{d}/R_{d}}\cong
\pi_{*}\oplus_{i\neq
j}\mathcal{H}om(\mathcal{Y}_{i},\mathcal{O}_{\mathbb{P}^{1}}/\mathcal{Y}_{j}).\]
Let us suppose that $\mathcal{E}\cong
\mathcal{O}_{\mathcal{Z}_{d}}\oplus 0\oplus \mathcal{O}\oplus
\mathcal{O}$, it is enough to consider one component of fixed
points by the symmetry of the computations, therefore,
\[0\rightarrow \overbrace{\mathcal{O}(-1,-d)}^{w_{0}}\oplus
\overbrace{\mathcal{O}}^{w_{1}}\oplus \overbrace{ 0}^{w_{2}}\oplus
\overbrace{0}^{w_{3}}\rightarrow \mathcal{O}^{4}\rightarrow
\mathcal{O}_{\mathcal{Z}_{d}}\oplus 0\oplus \mathcal{O}\oplus
\mathcal{O}\rightarrow 0, \]
\[\mbox{and}\ \mathcal{T}_{R_{d}}\cong
\pi_{*}\mathcal{H}om(\mathcal{O}(-1,-d)\oplus
\mathcal{O},\mathcal{O}_{\mathcal{Z}_{d}}\oplus 0\oplus
\mathcal{O}\oplus \mathcal{O} ).\] Definitely, what we have is
\begin{equation}\label{normal}\mathcal{N}_{\mathbb{P}^{d}/R_{d}}\cong
\pi_{*}\mathcal{O}_{\mathcal{Z}_{d}}\otimes
\mathcal{O}_{\mathbb{P}^{d}\times \mathbb{P}^{1}}(1,d)\oplus
(\pi_{*}\mathcal{O}_{\mathbb{P}^{d}\times
\mathbb{P}^{1}}(1,d))^{2}\oplus
(\pi_{*}\mathcal{O}_{\mathbb{P}^{d}\times \mathbb{P}^{1}})^{2}
.\end{equation} Since $\pi_{*}\mathcal{O}_{\mathcal{Z}_{d}}$ is a
bundle of rank $d$, the rank of the normal bundle is $3d+4$.
Therefore we need to compute,
$c_{3d+4}^{T}(\mathcal{N}_{\mathbb{P}^{d}/R_{d}})$ in the
equivariant Chow ring of $\mathbb{P}^{d}$, (see $\S$1.8 of
\cite{cf}),
\[A^{T}_{*}(\mathbb{P}^{d})=\mathbb{Z}(h,t)/\prod_{i=0}^{d}(h+w_{i}\,t).\]
It will be a polynomial of degree $3d+4$ in the variable $h$. The
fiber of the normal bundle is isomorphic to

\[\hom(\mathcal{O}_{\mathbb{P}^{1}},\mathcal{O}_{Z_{d}})\oplus \hom(\mathcal{O}_
{\mathbb{P}^{1}},\mathcal{O}_{\mathbb{P}^{1}})\oplus
\hom(\mathcal{O}_
{\mathbb{P}^{1}},\mathcal{O}_{\mathbb{P}^{1}})\oplus
\]\[\hom(\mathcal{O}_
{\mathbb{P}^{1}}(-d),\mathcal{O}_{\mathbb{P}^{1}})\oplus
\hom(\mathcal{O}_
{\mathbb{P}^{1}}(-d),\mathcal{O}_{\mathbb{P}^{1}}).\] and the
weights are by \cite{Bif}:
\[\underbrace{\hom(\mathcal{O}_{\mathbb{P}^{1}},\mathcal{O}_{Z_{d}})}_{w_{0}-w_{1}}
\oplus\underbrace{ \hom(\mathcal{O}_
{\mathbb{P}^{1}},\mathcal{O}_{\mathbb{P}^{1}})}_{w_{2}-w_{1}}\oplus
\underbrace{\hom(\mathcal{O}_
{\mathbb{P}^{1}},\mathcal{O}_{\mathbb{P}^{1}})}_{w_{3}-w_{1}}\oplus
\]\[\underbrace{\hom(\mathcal{O}_
{\mathbb{P}^{1}}(-d),\mathcal{O}_{\mathbb{P}^{1}})}_{w_{2}-w_{0}}\oplus
\underbrace{\hom(\mathcal{O}_
{\mathbb{P}^{1}}(-d),\mathcal{O}_{\mathbb{P}^{1}})}_{w_{3}-w_{0}}.\]

Given a $T-$equivariant vector bundle $E\rightarrow X$ we get a
canonical decomposition $E=\bigoplus_{\chi\in \hat{T}}E^{\chi}$
where $E^{\chi}$ denotes the eigensubbundle consisting of vectors in $E$ on which
$T$ acts with the character $\chi$.

The $i-$esima equivariant Chern class of a $T-$equivariant bundle
of rank $r$ over $\mathbb{P}^{d}$ is such that the action of $T$
over each fiber is given by the character. It is by $\S 2.2.1$,
\cite{cf}:
\[c_{i}^{T}=\sum_{j=0}^{i}{r-j \choose
i-j}c_{j}(E_{\chi_{i-j}})\chi_{i-j}.\] This formula relates the
equivariant Chern classes of a bundle with the usual Chern
classes. Thus, the problem is reduced to compute the usual Chern
classes of the normal bundle and by using the Whitney formula,
\[c_{3d+4}^{T}(\mathcal{N}_{\mathbb{P}^{d}/R_{d}})=c_{d}^{T}
(\pi_{*}\mathcal{O}_{\mathcal{Z}_{d}})\cdot\,c_{d+1}^{T}(\pi_{*}\mathcal{O}(1,d))\cdot\,
c_{d+1}^{T}(\pi_{*}\mathcal{O}(1,d))\,\cdot
c_{1}^{T}(\pi_{*}\mathcal{O})\cdot c_{1}^{T}(\pi_{*}\mathcal{O}).
\]
For computing the Chern classes of the equivariant subbundle
$\pi_{*}\mathcal{O}_{\mathcal{Z}_{d}}$ is required a little more
work. For this purpose, we consider the exact sequence
\[0\rightarrow \mathcal{O}_{\mathbb{P}^{d}\times \mathbb{P}^{1}}(-1,-d)\rightarrow
\mathcal{O}\rightarrow \mathcal{O}_{\mathcal{Z}_{d}}\rightarrow 0\
\ \ \rm{on} \ \mathbb{P}^{d}\times \mathbb{P}^{1},\]
\[0\rightarrow \pi_{*}\mathcal{O}_{\mathbb{P}^{d}\times \mathbb{P}^{1}}(-1,-d)\rightarrow
\pi_{*}\mathcal{O}_{\mathbb{P}^{d}\times
\mathbb{P}^{1}}\rightarrow
\pi_{*}\mathcal{O}_{\mathcal{Z}_{d}}\rightarrow
R^{1}\pi_{*}\mathcal{O}_{\mathbb{P}^{d}\times
\mathbb{P}^{1}}(-1,-d)\rightarrow 0.\] It follows
$\pi_{*}\mathcal{O}_{\mathbb{P}^{d}\times
\mathbb{P}^{1}}(-1,-d)=0$, since the fibers of this bundle are
isomorphic to $H^{0}(\mathcal{O}_{ \mathbb{P}^{1}}(-d))$ which are
0-dimensional vectorial spaces.
\[R^{1}\pi_{*}\mathcal{O}(-1,-d)=R^{1}\pi_{*}(\pi_{2}^{*}
\mathcal{O}_{\mathbb{P}^{1}}(-d)\otimes
\pi^{*}\mathcal{O}_{\mathbb{P}^{d}}(-1))=R^{1}\pi_{*}\pi_{2}^{*}
(\mathcal{O}_{\mathbb{P}^{1}}(-d))\otimes
\mathcal{O}_{\mathbb{P}^{3}}(-1)\] and $R^{1}\pi_{*}\pi_{2}^{*}
(\mathcal{O}_{\mathbb{P}^{1}}(-d))\cong \mathcal{O}^{d-1}$ by
Serre duality. Definitely, we have
\[0\rightarrow \mathcal{O}\rightarrow \pi_{*}\mathcal{O}_{\mathcal{Z}_{d}}\rightarrow
\mathcal{O}^{d-1}\otimes
\pi^{*}\mathcal{O}_{\mathbb{P}^{d}}(-1)\rightarrow 0\ \ \rm{on}\ \
\mathbb{P}^{d}\] and therefore, by applying the Whitney formula,
the total Chern class of the bundle is
\begin{center}
\fbox{$c_{t}(\pi_{*}\mathcal{O}_{\mathcal{Z}_{d}})=\prod_{i=0}^{d-1}(1-t)$}
\end{center}

\noindent Now we can compute the Chern equivariant class of each
equivariant subbundle of (\ref{normal}). Let $E_{w_{i}-w_{j}}$
denote the eigensubbundle consisting of vectors in
$\mathcal{N}_{\mathbb{P}^{d}/R_{d}}$ on which $\mathbb{C}^{*}$
acts with weight $w_{i}-w_{j}$.

\begin{equation}\label{chern-d}c_{d}^{T}(\pi_{*}\mathcal{O}_{\mathcal{Z}_{d}})=(w_{0}-w_{1})\,(h-(w_{0}-w_{1}))^{d-1}
\end{equation}
\[c_{d+1}^{T}(E_{w_{2}-w_{0}})=(h+w_{2}-w_{0})^{d+1},\]


\[c_{d+1}^{T}(E_{w_{3}-w_{0}})=(h+w_{3}-w_{0})^{d+1},\]

\[c_{1}^{T}(E_{w_{2}-w_{1}})=(w_{2}-w_{1}),\]
\[c_{1}^{T}(E_{w_{3}-w_{1}})=(w_{3}-w_{1}).\]
We have,
\[c_{3d+4}^{T}(\mathcal{N}_{\mathbb{P}^{d}/R_{d}})=(h+(w_{2}-w_{0}))^{d+1}\cdot
 (h+(w_{3}-w_{0}))^{d+1}\cdot(w_{0}-w_{1})\cdot (h-(w_{0}-w_{1}))^{d-1}\]\[(w_{2}-w_{1})\cdot(w_{3}-w_{1}).\]

\vspace{0.5cm} {\bf We now study
$c_{d_{F}}(\mathcal{N}_{F/R_{d}})$, for the components of fixed
points of the second kind}.

We consider the varieties of fixed points isomorphic to
$\mathbb{P}^{b}\times \mathbb{P}^{a}$, with $b+a=d$, $b\geq a>0$.

\noindent Let be the component of fixed points defined by the
universal quotient in $\mathbb{P}^{b}\times \mathbb{P}^{a}\times
\mathbb{P}^{1}$:
\[0\rightarrow \mathcal{O}
(-1,0,-b)\oplus \mathcal{O}(0,-1,-a)\oplus 0\oplus 0\rightarrow
\mathcal{O}^{4}\rightarrow \mathcal{O}_{\mathcal{Z}_{b}}\oplus
\mathcal{O}_{\mathcal{Z}_{a}}\oplus \mathcal{O}\oplus
\mathcal{O}\rightarrow 0  \] \[\mathcal{N}_{\mathbb{P}^{b}\times
\mathbb{P}^{a}/R_{d}}=\pi_{12*}\mathcal{H}om(\mathcal{O}(-1,0,-b),\mathcal{O}_{\mathcal{Z}_{a}})
\oplus
\pi_{12*}\mathcal{H}om(\mathcal{O}(0,-1,-a),\mathcal{O}_{\mathcal{Z}_{b}})\oplus\]
\[\pi_{12*}\mathcal{H}om(\mathcal{O}(-1,0,-b),\mathcal{O})^{2}\oplus
\pi_{12*}\mathcal{H}om(\mathcal{O}(0,-1,-a),\mathcal{O})^{2}.\] We
denote $E_{w_{2}-w_{0}}$ the subbundle
$\pi_{12_{*}}(\mathcal{O}_{\mathbb{P}^{b}\times
\mathbb{P}^{a}\times \mathbb{P}^{1}}(1,0,b))$ of
$\mathcal{N}_{\mathbb{P}^{b}\times \mathbb{P}^{a}/R_{d}}$ on which
$\mathbb{C}^{*}$ acts with weight $w_{2}-w_{0}$. We compute its
equivariant Chern classes:
\[c_{b+1}^{T}(\underbrace{\pi_{12*}(\mathcal{O}_{\mathbb{P}^{b}\times \mathbb{P}^{a}\times \mathbb{P}^{1}}(1,0,b))}_{E_{w_{2}-w_{0}}}) =
(H+(w_{2}-w_{0}))^{b+1},\]
\[c_{b+1}^{T}(\underbrace{\pi_{12*}(\mathcal{O}_{\mathbb{P}^{b}\times \mathbb{P}^{a}\times \mathbb{P}^{1}}(1,0,b))}_{E_{w_{3}-w_{0}}})=
(H+(w_{3}-w_{0}))^{b+1},\]
\[c_{a+1}^{T}(\pi_{12*}(\underbrace{\mathcal{O}_{\mathbb{P}^{b}\times \mathbb{P}^{a}\times \mathbb{P}^{1}}(0,1,a))}_{E_{w_{2}-w_{0}}})=
(h+(w_{2}-w_{0}))^{a+1},\]
\[c_{a+1}^{T}(\pi_{12*}(\underbrace{\mathcal{O}_{\mathbb{P}^{b}\times \mathbb{P}^{a}\times \mathbb{P}^{1}}(0,1,a))}_{E_{w_{3}-w_{1}}})=
(h+(w_{3}-w_{1}))^{a+1},\] $\mathbb{C}^{*}$ acts on
$\pi_{12*}(\mathcal{O}_{\mathcal{Z}_{b}}\otimes
\mathcal{O}(0,1,a))$ with weight $w_{0}-w_{1}$. For computing its
equivariant Chern class $c^{T}_{b}(E_{w_{0}-w_{1}})$ we consider
the exact sequence,
\[0\rightarrow \mathcal{O}_{\mathbb{P}^{b}\times \mathbb{P}^{a}\times \mathbb{P}^{1}}(-1,0,-b)
\otimes\mathcal{O}(0,1,a)\rightarrow
\mathcal{O}_{\mathbb{P}^{b}\times \mathbb{P}^{a}\times
\mathbb{P}^{1}}(0,1,a) \rightarrow
\mathcal{O}_{\mathcal{Z}_{b}}\otimes \mathcal{O}(0,1,a)\rightarrow
0
\]
\ \begin{equation} \label{ab} 0\rightarrow
\pi_{12*}(\mathcal{O}_{\mathbb{P}^{b}\times \mathbb{P}^{a}\times
\mathbb{P}^{1}}(-1,1,-b+a))\rightarrow
\pi_{12*}(\mathcal{O}_{\mathbb{P}^{b}\times \mathbb{P}^{a}\times
\mathbb{P}^{1}}(0,1,a))\rightarrow \end{equation}

$$
\label{ab}\rightarrow
\pi_{12*}(\mathcal{O}_{\mathcal{Z}_{b}}\otimes
\mathcal{O}_{\mathbb{P}^{b}\times \mathbb{P}^{a}\times
\mathbb{P}^{1}}(0,1,a))\rightarrow
R^{1}\pi_{12*}\mathcal{O}_{\mathbb{P}^{b}\times
\mathbb{P}^{a}\times \mathbb{P}^{1}}(-1,1,-b+a))\rightarrow 0 $$
\[\pi_{12*}(\mathcal{O}_{\mathbb{P}^{b}\times \mathbb{P}^{a}\times
\mathbb{P}^{1}}(0,1,a))\cong \mathcal{O}(0,1)^{a+1}\] We suppose
$b>a$, thus we have

\noindent $\pi_{12*}\mathcal{O}_{\mathbb{P}^{b}\times \mathbb{P}^{a}\times
\mathbb{P}^{1}}(-1,1,-b+a)=0$, 
$R^{1}\pi_{12*}\mathcal{O}_{\mathbb{P}^{b}\times
\mathbb{P}^{a}\times \mathbb{P}^{1}}(-1,1,-b+a)\neq 0$ and
$$R^{1}\pi_{12*}\mathcal{O}_{\mathbb{P}^{b}\times \mathbb{P}^{a}\times
\mathbb{P}^{1}}(-1,1,-b+a)=R^{1}\pi_{12*}(\pi^{*}_{12}\mathcal{O}(-1,1)\otimes
\pi_{3}^{*}\mathcal{O}(-b+a))=$$

\noindent $\mathcal{O}_{\mathbb{P}^{b}\times
\mathbb{P}^{a}}(-1,1)\otimes
R^{1}\pi_{12*}\pi_{3}^{*}\mathcal{O}(-b+a)\cong
\mathcal{O}_{\mathbb{P}^{b}\times \mathbb{P}^{a}}(-1,1)\otimes
\pi_{12*}\pi_{3}^{*}\mathcal{O}_{\mathbb{P}^{1}}(b-a-2)=\mathcal{O}_{\mathbb{P}^{b}\times
\mathbb{P}^{a}}(-1,1)\otimes \mathcal{O}^{b-a-1}$.

\noindent In the case $a=b$, we have that
\begin{equation}0\rightarrow \mathcal{O}_{\mathbb{P}^{b}\times \mathbb{P}^{a}}(0,1)\rightarrow
\pi_{12*}\mathcal{O}_{\mathbb{P}^{b}\times\mathbb{P}^{a}\times
\mathbb{P}^{1}}(0,1,a) \rightarrow
\pi_{12*}(\mathcal{O}_{\mathcal{Z}_{b}}\otimes
\mathcal{O}_{\mathbb{P}^{b}\times \mathbb{P}^{a}\times
\mathbb{P}^{1}}(0,1,a))\rightarrow 0.\end{equation}

\noindent In either case, by applying Whitney formula to
(\ref{ab}), we see that,
$$c^{T}(E_{w_{0}-w_{1}})=c_{b}^{T}(\pi_{12*}(\mathcal{O}_{\mathcal{Z}_{b}}\otimes
\mathcal{O}(0,1,a)))= (-H+h+w_{0}-w_{1})^{b-a-1}\cdot
(h+w_{0}-w_{1})^{a+1}.$$



\noindent $\mathbb{C}^{*}$ acts on
$\pi_{12*}(\mathcal{O}_{\mathcal{Z}_{a}}\otimes
\mathcal{O}_{\mathbb{P}^{b}\times \mathbb{P}^{a}\times
\mathbb{P}^{1}}(1,0,b))$ with weight $w_{1}-w_{0}$. We now compute
$c_{a}^{T}(\pi_{12*}(\mathcal{O}_{\mathcal{Z}_{a}}\otimes
\mathcal{O}(1,0,b)))$. We consider also the case in which $a=b$.
Consider the exact sequence
\[0\rightarrow \mathcal{O}_{\mathbb{P}^{b}\times \mathbb{P}^{a}\times
\mathbb{P}^{1}}(0,-1,-a)\otimes\mathcal{O}(1,0,b)\rightarrow
\mathcal{O}(1,0,b) \rightarrow
\mathcal{O}_{\mathcal{Z}_{a}}\otimes \mathcal{O}(1,0,b)\rightarrow
0,
\]
\[0\rightarrow
\pi_{12*}(\mathcal{O}_{\mathbb{P}^{b}\times \mathbb{P}^{a}\times
\mathbb{P}^{1}}(1,-1,b-a))\rightarrow
\pi_{12*}(\mathcal{O}(1,0,b))\rightarrow
\pi_{12*}(\mathcal{O}_{\mathcal{Z}_{a}}\otimes
\mathcal{O}(1,0,b))\rightarrow\]\begin{equation}
\label{ba}\rightarrow
R^{1}\pi_{12*}\mathcal{O}_{\mathbb{P}^{b}\times
\mathbb{P}^{a}\times \mathbb{P}^{1}}(1,-1,b-a))\rightarrow
0.\end{equation} We observe in this case
\[R^{1}\pi_{12*}\mathcal{O}_{\mathbb{P}^{b}\times \mathbb{P}^{a}\times
\mathbb{P}^{1}}(1,-1,b-a)\cong \mathcal{O}(1,-1)\otimes
\pi_{12*}(\pi^{*}_{3}\mathcal{O}_{\mathbb{P}^{1}}(a-b-2))=0,\]
\begin{equation}
\pi_{12*}(\mathcal{O}_{\mathbb{P}^{b}\times \mathbb{P}^{a}\times
\mathbb{P}^{1}}(1,-1,b-a))\cong \mathcal{O}(1,-1)\otimes
\mathcal{O}^{b-a+1},
\end{equation}
\begin{equation}
\pi_{12*}(\mathcal{O}_{\mathbb{P}^{b}\times \mathbb{P}^{a}\times
\mathbb{P}^{1}}(1,0,b))\cong \mathcal{O}(1,0)\otimes
\mathcal{O}^{b+1},
\end{equation}
$\pi_{12*}(\mathcal{O}_{\mathbb{P}^{b}\times \mathbb{P}^{a}\times
\mathbb{P}^{1}}(1,-1,b-a))$ is a bundle with total Chern class:
\[c_{t}(\pi_{12*}\mathcal{O}_{\mathbb{P}^{b}\times \mathbb{P}^{a}\times
\mathbb{P}^{1}}(1,-1,b-a))=((H-h)t+1)^{b-a+1},\]
\[c_{t}^{T}(\pi_{12*}\mathcal{O}_{\mathbb{P}^{b}\times \mathbb{P}^{a}\times
\mathbb{P}^{1}}(1,-1,b-a))=((H-h)t+w_{1}-w_{0})^{b-a+1},\]
\[c_{t}(\pi_{12*}(\mathcal{O}_{\mathbb{P}^{b}\times \mathbb{P}^{a}\times
\mathbb{P}^{1}}(1,0,b)))=(Ht+1)^{b+1},\]
\[c_{t}^{T}(\pi_{12*}(\mathcal{O}_{\mathbb{P}^{b}\times \mathbb{P}^{a}\times
\mathbb{P}^{1}}(1,0,b)))=(Ht+w_{1}-w_{0})^{b+1},\] By applying
Whitney formula to (\ref{ba}), we get that
\[c_{a}(E_{w_{1}-w_{0}})=\frac{(H+w_{1}-w_{0})^{b+1}}{(H-h+w_{1}-w_{0})^{b-a+1}},\]
therefore,

\[c_{3d+4}^{T}(\mathcal{N}_{\mathbb{P}^{b}\times
\mathbb{P}^{a}}/R_{d})=c^{T}_{b+1}(E_{w_{2}-w_{0}})\cdot
c^{T}_{b+1}(E_{w_{3}-w_{0}})\cdot
c^{T}_{a+1}(E_{w_{2}-w_{0}})\cdot
c^{T}_{a+1}(E_{w_{3}-w_{1}})\cdot\]\[
c_{b}^{T}(E_{w_{0}-w_{1}})\cdot c_{a}^{T}(E_{w_{1}-w_{0}}),\]

\noindent that is, for $b-a\geq 1$,
$$c_{3d+4}^{T}(\mathcal{N}_{\mathbb{P}^{b}\times
\mathbb{P}^{a}}/R_{d})=(H+w_{2}-w_{0})^{b+1}\cdot
(H+w_{3}-w_{0})^{b+1}\cdot (h+w_{2}-w_{0})^{a+1}\cdot$$$$
(-H+h+w_{0}-w_{1})^{b-a-1}\cdot (h+w_{0}-w_{1})^{a+1}\cdot
(H+w_{1}-w_{0})^{b+1}\cdot (H-h+w_{1}-w_{0})^{a-b-1},$$ and for
$a=b$,
$$c_{3d+4}^{T}(\mathcal{N}_{\mathbb{P}^{a}\times
\mathbb{P}^{a}}/R_{d})=(H+w_{2}-w_{0})^{a+1}\cdot
(H+w_{3}-w_{0})^{a+1}\cdot
(h+w_{2}-w_{0})^{a+1}\cdot(h+w_{0}-w_{1})^{a+1}\cdot$$$$
(H+w_{1}-w_{0})^{a+1}\cdot (H-h+w_{1}-w_{0})^{-1}.$$

\cqd
\section{Appendix A: calculation of Pl\"ucker degree of  $R_{3}$.}


We want to compute the degree of $R_{3}$ by the morphism induced
by the divisor $\alpha$, that is, the generalized Pl\"ucker
embbeding \cite{RRW}. The intersection we compute, is
\[P_{3}=\int_{R_{3}}(\alpha^{16}\cap [R_{3}]).\]

We apply Bott residue formula. We have 24 summands, one for each
component of fixed points. We know what the denominator is
(\ref{fnormal}), and the restrictions of $\alpha$ to each
subvariety of fixed points. The 24 summands corresponding to the
24 components of fixed points are:

\begin{enumerate}
\item \[Hilb^{0}_{\mathbb{P}^{1}}\times
Hilb^{3}_{\mathbb{P}^{1}}\times Hilb^{t+1}_{\mathbb{P}^{1}}\times
Hilb^{t+1}_{\mathbb{P}^{1}}\]
\[\frac{(h+w_{2}+w_{3})^{16}}{(h+w_{2}-w_{1})^{4}
 (h+w_{3}-w_{1})^{4}(w_{1}-w_{0})(h-w_{1}+w_{0})^{2}(w_{2}-w_{0})(w_{3}-w_{0})}\]
\item \[Hilb^{t+1}_{\mathbb{P}^{1}}\times
Hilb^{3}_{\mathbb{P}^{1}}\times Hilb^{t+1}_{\mathbb{P}^{1}}\times
Hilb^{0}_{\mathbb{P}^{1}}\]
\[\frac{(h+w_{0}+w_{2})^{16}}{(h+w_{0}-w_{1})^{4}
 (h+w_{2}-w_{1})^{4}(w_{1}-w_{3})(h-w_{1}+w_{3})^{2}(w_{2}-w_{3})(w_{0}-w_{3})}\]

\item \[Hilb^{t+1}_{\mathbb{P}^{1}}\times
Hilb^{3}_{\mathbb{P}^{1}}\times Hilb^{0}_{\mathbb{P}^{1}}\times
Hilb^{t+1}_{\mathbb{P}^{1}}\]
\[\frac{(h+w_{0}+w_{3})^{16}}{(h+w_{0}-w_{1})^{4}
 (h+w_{3}-w_{1})^{4}(w_{1}-w_{2})(h-w_{1}+w_{2})^{2}(w_{0}-w_{2})(w_{3}-w_{2})}\]

 \item \[Hilb^{0}_{\mathbb{P}^{1}}\times Hilb^{t+1}_{\mathbb{P}^{1}}\times
Hilb^{3}_{\mathbb{P}^{1}}\times Hilb^{t+1}_{\mathbb{P}^{1}}\]
\[\frac{(h+w_{1}+w_{3})^{16}}{(h+w_{1}-w_{2})^{4}
 (h+w_{3}-w_{2})^{4}(w_{2}-w_{0})(h-w_{2}+w_{0})^{2}(w_{1}-w_{0})(w_{3}-w_{0})}\]

 \item \[Hilb^{t+1}_{\mathbb{P}^{1}}\times Hilb^{t+1}_{\mathbb{P}^{1}}\times
Hilb^{3}_{\mathbb{P}^{1}}\times Hilb^{0}_{\mathbb{P}^{1}}\]
\[\frac{(h+w_{0}+w_{1})^{16}}{(h+w_{0}-w_{2})^{4}
 (h+w_{1}-w_{2})^{4}(w_{2}-w_{3})(h-w_{2}+w_{3})^{2}(w_{0}-w_{3})(w_{1}-w_{3})}\]

\item \[Hilb^{t+1}_{\mathbb{P}^{1}}\times
Hilb^{0}_{\mathbb{P}^{1}}\times Hilb^{3}_{\mathbb{P}^{1}}\times
Hilb^{t+1}_{\mathbb{P}^{1}}\]
\[\frac{(h+w_{0}+w_{3})^{16}}{(h+w_{0}-w_{2})^{4}
 (h+w_{3}-w_{2})^{4}(w_{2}-w_{1})(h-w_{2}+w_{1})^{2}(w_{0}-w_{1})(w_{3}-w_{1})}\]


\item \[Hilb^{3}_{\mathbb{P}^{1}}\times
Hilb^{t+1}_{\mathbb{P}^{1}}\times Hilb^{0}_{\mathbb{P}^{1}}\times
Hilb^{t+1}_{\mathbb{P}^{1}}\]
\[\frac{(h+w_{1}+w_{2})^{16}}{(h+w_{1}-w_{0})^{4}
 (h+w_{2}-w_{0})^{4}(w_{0}-w_{3})(h-w_{0}+w_{3})^{2}(w_{1}-w_{3})(w_{2}-w_{3})}\]

\item \[Hilb^{3}_{\mathbb{P}^{1}}\times
Hilb^{0}_{\mathbb{P}^{1}}\times Hilb^{t+1}_{\mathbb{P}^{1}}\times
Hilb^{t+1}_{\mathbb{P}^{1}}\]
\[\frac{(h+w_{2}+w_{3})^{16}}{(h+w_{2}-w_{0})^{4}
 (h+w_{3}-w_{0})^{4}(w_{0}-w_{1})(h-w_{0}+w_{1})^{2}(w_{2}-w_{1})(w_{3}-w_{1})}\]

\item \[Hilb^{3}_{\mathbb{P}^{1}}\times
Hilb^{t+1}_{\mathbb{P}^{1}}\times Hilb^{0}_{\mathbb{P}^{1}}\times
Hilb^{t+1}_{\mathbb{P}^{1}}\]
\[\frac{(h+w_{1}+w_{3})^{16}}{(h+w_{1}-w_{0})^{4}
 (h+w_{3}-w_{0})^{4}(w_{0}-w_{2})(h-w_{0}+w_{2})^{2}(w_{1}-w_{2})(w_{3}-w_{2})}\]

\item \[Hilb^{0}_{\mathbb{P}^{1}}\times
Hilb^{t+1}_{\mathbb{P}^{1}}\times
Hilb^{t+1}_{\mathbb{P}^{1}}\times Hilb^{3}_{\mathbb{P}^{1}}\]
\[\frac{(h+w_{1}+w_{2})^{16}}{(h+w_{1}-w_{3})^{4}
 (h+w_{2}-w_{3})^{4}(w_{3}-w_{0})(h-w_{3}+w_{0})^{2}(w_{1}-w_{0})(w_{2}-w_{0})}\]

\item \[Hilb^{t+1}_{\mathbb{P}^{1}}\times
Hilb^{0}_{\mathbb{P}^{1}}\times Hilb^{t+1}_{\mathbb{P}^{1}}\times
Hilb^{3}_{\mathbb{P}^{1}}\]
\[\frac{(h+w_{0}+w_{2})^{16}}{(h+w_{0}-w_{3})^{4}
 (h+w_{2}-w_{3})^{4}(w_{3}-w_{1})(h-w_{3}+w_{1})^{2}(w_{0}-w_{1})(w_{2}-w_{1})}\]

\item \[Hilb^{t+1}_{\mathbb{P}^{1}}\times
Hilb^{t+1}_{\mathbb{P}^{1}}\times Hilb^{0}_{\mathbb{P}^{1}}\times
Hilb^{3}_{\mathbb{P}^{1}}\]

\[\frac{(h+w_{0}+w_{1})^{16}}{(h+w_{0}-w_{3})^{4}
 (h+w_{1}-w_{3})^{4}(w_{3}-w_{2})(h-w_{3}+w_{2})^{2}(w_{0}-w_{2})(w_{1}-w_{2})}\]

\item \[Hilb^{1}_{\mathbb{P}^{1}}\times
Hilb^{2}_{\mathbb{P}^{1}}\times Hilb^{t+1}_{\mathbb{P}^{1}}\times
Hilb^{t+1}_{\mathbb{P}^{1}}\]
\scriptsize
\[\frac{(H+h+w_{2}+w_{3})^{16}}{(h+w_{1}-w_{0})^{2}
 (H+2h+w_{0}-w_{1})(h+w_{2}-w_{0})^{2}(h+w_{3}-w_{0})^{2}(H+w_{3}-w_{1})^{3}(H+w_{2}-w_{1})^{3}}\]

\normalsize


\item \[Hilb^{1}_{\mathbb{P}^{1}}\times
Hilb^{t+1}_{\mathbb{P}^{1}}\times Hilb^{2}_{\mathbb{P}^{1}}\times
Hilb^{t+1}_{\mathbb{P}^{1}}\]
\scriptsize \[\frac{(H+h+w_{1}+w_{3})^{16}}{(h+w_{2}-w_{0})^{2}
 (H+2h+w_{0}-w_{2})(h+w_{1}-w_{0})^{2}(h+w_{3}-w_{0})^{2}(H+w_{1}-w_{2})^{3}(H+w_{3}-w_{2})^{3}}\]
\normalsize

\item \[Hilb^{1}_{\mathbb{P}^{1}}\times
Hilb^{t+1}_{\mathbb{P}^{1}}\times Hilb^{2}_{\mathbb{P}^{1}}\times
Hilb^{t+1}_{\mathbb{P}^{1}}\]
\scriptsize \[\frac{(H+h+w_{1}+w_{2})^{16}}{(h+w_{3}-w_{0})^{2}
 (H+2h+w_{0}-w_{3})(h+w_{2}-w_{0})^{2}(h+w_{1}-w_{0})^{2}(H+w_{2}-w_{3})^{3}(H+w_{1}-w_{0})^{3}}\]
\normalsize

\item \[Hilb^{t+1}_{\mathbb{P}^{1}}\times
Hilb^{1}_{\mathbb{P}^{1}}\times Hilb^{2}_{\mathbb{P}^{1}}\times
Hilb^{t+1}_{\mathbb{P}^{1}}\]
\scriptsize \[\frac{(H+h+w_{0}+w_{3})^{16}}{(h+w_{2}-w_{1})^{2}
 (H+2h+w_{1}-w_{2})(h+w_{0}-w_{1})^{2}(h+w_{3}-w_{1})^{2}(H+w_{0}-w_{2})^{3}(H+w_{3}-w_{2})^{3}}\]
\normalsize


\item \[Hilb^{t+1}_{\mathbb{P}^{1}}\times
Hilb^{1}_{\mathbb{P}^{1}}\times Hilb^{t+1}_{\mathbb{P}^{1}}\times
Hilb^{2}_{\mathbb{P}^{1}}\]
\scriptsize \[\frac{(H+h+w_{0}+w_{2})^{16}}{(h+w_{3}-w_{1})^{2}
 (H+2h+w_{1}-w_{3})(h+w_{0}-w_{1})^{2}(h+w_{2}-w_{1})^{2}(H+w_{0}-w_{3})^{3}(H+w_{2}-w_{3})^{3}}\]
\normalsize


\item \[Hilb^{2}_{\mathbb{P}^{1}}\times
Hilb^{1}_{\mathbb{P}^{1}}\times Hilb^{t+1}_{\mathbb{P}^{1}}\times
Hilb^{t+1}_{\mathbb{P}^{1}}\]

\scriptsize \[\frac{(H+h+w_{2}+w_{3})^{16}}{(h+w_{3}-w_{1})^{2}
 (H+2h+w_{1}-w_{3})(h+w_{2}-w_{1})^{2}(h+w_{3}-w_{1})^{2}(H+w_{2}-w_{0})^{3}(H+w_{3}-w_{0})^{3}}\]
\normalsize


\item \[Hilb^{t+1}_{\mathbb{P}^{1}}\times
Hilb^{2}_{\mathbb{P}^{1}}\times Hilb^{1}_{\mathbb{P}^{1}}\times
Hilb^{t+1}_{\mathbb{P}^{1}}\]

\scriptsize \[\frac{(H+h+w_{0}+w_{3})^{16}}{(h+w_{0}-w_{2})^{2}
 (H+2h+w_{2}-w_{0})(h+w_{0}-w_{2})^{2}(h+w_{3}-w_{2})^{2}(H+w_{0}-w_{1})^{3}(H+w_{3}-w_{1})^{3}}\]
\normalsize

\item \[Hilb^{t+1}_{\mathbb{P}^{1}}\times
Hilb^{t+1}_{\mathbb{P}^{1}}\times Hilb^{1}_{\mathbb{P}^{1}}\times
Hilb^{2}_{\mathbb{P}^{1}}\]

\scriptsize \[\frac{(H+h+w_{0}+w_{3})^{16}}{(h+w_{0}-w_{2})^{2}
 (H+2h+w_{2}-w_{0})(h+w_{0}-w_{2})^{2}(h+w_{3}-w_{2})^{2}(H+w_{0}-w_{1})^{3}(H+w_{3}-w_{1})^{3}}\]
\normalsize


\item \[Hilb^{2}_{\mathbb{P}^{1}}\times
Hilb^{t+1}_{\mathbb{P}^{1}}\times Hilb^{1}_{\mathbb{P}^{1}}\times
Hilb^{t+1}_{\mathbb{P}^{1}}\] \scriptsize
\[\frac{(H+h+w_{0}+w_{1})^{16}}{(h+w_{0}-w_{2})^{2}
 (H+2h+w_{2}-w_{0})(h+w_{1}-w_{2})^{2}(h+w_{3}-w_{2})^{2}(H+w_{1}-w_{0})^{3}(H+w_{3}-w_{0})^{3}}\]
\normalsize


\item \[Hilb^{t+1}_{\mathbb{P}^{1}}\times
Hilb^{t+1}_{\mathbb{P}^{1}}\times Hilb^{2}_{\mathbb{P}^{1}}\times
Hilb^{1}_{\mathbb{P}^{1}}\]

\scriptsize \[\frac{(H+h+w_{0}+w_{1})^{16}}{(h+w_{2}-w_{1})^{2}
 (H+2h+w_{1}-w_{2})(h+w_{0}-w_{3})^{2}(h+w_{1}-w_{3})^{2}(H+w_{0}-w_{2})^{3}(H+w_{1}-w_{2})^{3}}\]
\normalsize

\item \[Hilb^{t+1}_{\mathbb{P}^{1}}\times
Hilb^{2}_{\mathbb{P}^{1}}\times Hilb^{t+1}_{\mathbb{P}^{1}}\times
Hilb^{1}_{\mathbb{P}^{1}}\]
\scriptsize \[\frac{(H+h+w_{0}+w_{2})^{16}}{(h+w_{1}-w_{3})^{2}
 (H+2h+w_{3}-w_{1})(h+w_{0}-w_{3})^{2}(h+w_{2}-w_{3})^{2}(H+w_{0}-w_{1})^{3}(H+w_{2}-w_{3})^{3}}\]
\normalsize

\item \[Hilb^{2}_{\mathbb{P}^{1}}\times
Hilb^{t+1}_{\mathbb{P}^{1}}\times
Hilb^{t+1}_{\mathbb{P}^{1}}\times Hilb^{1}_{\mathbb{P}^{1}}\]

\scriptsize \[\frac{(H+h+w_{1}+w_{2})^{16}}{(h+w_{0}-w_{3})^{2}
 (H+2h+w_{3}-w_{0})(h+w_{2}-w_{3})^{2}(h+w_{1}-w_{3})^{2}(H+w_{1}-w_{0})^{3}(H+w_{2}-w_{0})^{3}}\]
\normalsize

\end{enumerate}

Once we have all the summands, we take the direct image by the
morphism $\pi_{1}:\mathbb{P}^{3}\times \mathbb{P}^{1}\rightarrow
\mathbb{P}^{3}$ and $\pi_{12}:\mathbb{P}^{2}\times
\mathbb{P}^{1}\times \mathbb{P}^{1}\rightarrow
\mathbb{P}^{2}\times \mathbb{P}^{1}$ for the first kind of
components and for the second kind of components respectively. The
only terms surviving in the first case are those in $h^{3}$, and
in the second case the terms in $H^{2}h$. We have used Maple
program to make the computations.  This degree is 128. This result
coincides with the one obtained by means of the Vafa-Intriligator
formula, \cite{Ber1} and indeed it follows easily from Vafa-Intriligator
formula that the degree $P_{d}$ coincides with $2^{2d+1}$. 




\end{document}